%% file: surfim-arxiv.tex
\documentclass{article} 
\usepackage{graphicx}

\usepackage[greek,english]{babel}

\usepackage[utf8x]{inputenc}

\usepackage{amsmath,amsthm,amsfonts,amssymb,epsfig,color}
\usepackage{cite}
\usepackage{bm}  

\usepackage{enumerate}

\usepackage{hyperref}  
\hypersetup{colorlinks,%
            citecolor=green,%
            filecolor=black,%
            linkcolor=red,%
            urlcolor=blue}

\newtheorem{corollary}{Corollary}
\newtheorem{theorem}{Theorem}
\newtheorem{lemma}{Lemma}
\newtheorem{observation}{Observation}

\newcommand{\es}{^{o}}

\newcommand{\tw}{{\mathbf{tw}}}

\newcommand{\eg}{{\mathbf{eg}}}

\newcommand{\n}{{\mathbf{n}}}
\newcommand{\m}{{\mathbf{m}}}
\newcommand{\im}{\leq_{\text{im}}}
\newcommand{\img}{{\mathbf{Im}}}

\begin{document}

\title{Excluding Graphs as Immersions in Surface Embedded Graphs}

\date{}

\author{Archontia C. Giannopoulou\footnotemark[2]~\footnotemark[5]~~\footnotemark[6]~, Marcin Kami\'{n}ski\footnotemark[3]~,\\ and Dimitrios M. Thilikos\footnotemark[4]~\footnotemark[7]}

\renewcommand{\thefootnote}{\fnsymbol{footnote}}

\footnotetext[2]{Department of Informatics, University of Bergen, P.O. Box 7803, N-5020 Bergen, Norway. Email:
Archontia.Giannopoulou@ii.uib.no}
\footnotetext[3]{D\'{e}partement d'Informatique, Universit\'{e} Libre de Bruxelles
and Instytut Informatyki, Uniwersytet Warszawski.
Email: mjk@mimuw.edu.pl}
\footnotetext[4]{CNRS, LIRMM and Department of Mathematics, National and Kapodistrian University of Athens.
Email: sedthilk@thilikos.info}
\footnotetext[5]{Partially supported by European Research Council (ERC) Grant "Rigorous Theory
of Preprocessing", reference 267959.}
\footnotetext[6]{Partially supported by a grant of the Special Account for Research Grants of the  National and Kapodistrian University of Athens (project code: 70/4/10311).}
\footnotetext[7]{Co-financed by the European Union (European Social Fund - ESF) and Greek national funds through the Operational Program ``Education and Lifelong Learning'' of the National Strategic Reference Framework (NSRF) - Research Funding Program: ``Thales. Investing in knowledge society through the European Social Fund.''}

\renewcommand{\thefootnote}{\arabic{footnote}}

\maketitle

\begin{abstract}
We prove a structural characterization of graphs that forbid a fixed graph $H$
as an immersion and can be embedded in a surface of Eüler genus $\gamma$. 
In particular, we prove that a graph $G$ that excludes some connected graph $H$ as 
an immersion and is embedded in a surface of Eüler genus $\gamma$ has either ``small" 
treewidth (bounded by a function of $H$ and $\gamma$) or “small” edge connectivity 
(bounded by the maximum degree of $H$). Using the same techniques we also prove
an excluded grid theorem on bounded genus graphs for the immersion relation.
\end{abstract}

\noindent \textbf{Keywords:} Surface Embeddable Graphs, Immersion Relation, Treewidth, Edge Connectivity.

\section{Introduction}

A graph $H$ is an immersion of a graph $G$ if it can be obtained from $G$ by removing vertices or edges, and splitting off adjacent pairs of edges.
The class of all graphs was proved to be  well-quasi-ordered under the 
the immersion relation by Robertson and Seymour in the 
last paper of their Graph Minors series~\cite{RobertsonS10}. Certainly, this work was 
mostly dedicated to minors and not immersions and has been the source of  
many theorems regarding the structure of graphs excluding some graph $H$ as a minor. Moreover, the minor relation has been extensively studied
the past two decades and many structural results have been proven for minors with interesting algorithmic consequences (see, for example,~\cite{RobertsonS-XVI,RobertsonS-GMXIII,RobertsonST94,DemaineFHT05sube,Kostochka84,Thomason01}).
However, structural results for immersions started appearing only recently.
In 2011, DeVos et al. proved that if the minimum degree of a graph $G$ is $200t,$ then $G$ contains the complete graph on $t$ vertices as an immersion~\cite{2011arXiv1101.2630D}.
In~\cite{FerraraGTW08} Ferrara et al.,
provided a lower bound (depending on graph $H$) on the minimum degree of a
graph $G$ that ensures that $H$ is contained in $G$ as an immersion.
Furthermore, Wollan recently proved a structural theorem for graphs excluding complete graphs as immersions
as well as a sufficient condition such that any graph which satisfies the condition admits a wall as an immersion~\cite{abs-1302-3867}. The result in~\cite{abs-1302-3867} can be seen as  an immersion counterpart of the grid exclusion theorem~\cite{RobertsonST94}, stated for walls instead of grids and using an alternative graph parameter instead of treewidth.

In terms of graph colorings, Abu-Khzam and Langston in~\cite{Abu-KhzamL03} provided evidence supporting the immersion ordering analog of Hadwiger's Conjecture, that is, the conjecture stating that if the chromatic number of a graph $G$ is at least $t$, then $G$ contains the complete graph on $t$ vertices as an immersion, and proved it for $t\leq 4$. For $t=5,6,7$, see~\cite{Lescure1988325,1213.05137}. 
For algorithmic results on immersions, see~\cite{KawarabayashiK12,BelmonteHKPT12,GiannopoulouSZ12,GroheKMW11}.

In this paper, we prove structural results for the immersion relation on graphs embeddable  on a fixed surface.
In particular, we show that if $G$ is a graph that is embeddable on a surface of Eüler genus $\gamma$ and
$H$ is a connected graph then one of the following is true:
either $G$ has bounded treewidth (by a function that depends only on $\gamma$ and $H$), or its edge connectivity is bounded by the maximum degree of $H$, or it contains $H$ as a (strong) immersion. 
Furthermore, we refine our results 
to obtain a counterpart of the grid exclusion theorem for 
immersions. In particular, 
we prove (Theorem~\ref{thm:exclgrdimrs}) that there exists a function $f:\Bbb{N}\rightarrow\Bbb{N}$ such that if $G$ is a $4$-edge-connected graph embedded on a surface of Eüler genus $\gamma$ and the treewidth of $G$ is at least $f(\gamma)\cdot k$, then $G$ contains the $k\times k$-grid as an immersion.
Notice that  the edge connectivity requirement is necessary here as big treewidth alone is not enough to ensure the existence of a graph with a vertex of degree 4 as an immersion. Although a wall of height at least $h$ has treewidth at least $h$, it does not contain the complete graph on $t$ vertices as an immersion, for any $t\geq 5$.
Finally, our results imply that when restricted to graphs of sufficiently  big treewidth  embeddable on a fixed surface, large edge connectivity forces the existence of a large clique as an immersion.

Our result reveals several aspects of the behavior of the immersion relation on surface embeddable 
graphs.  The proofs  exploit variants of the grid exclusion theorem for surfaces 
proved in~\cite{FominGT11} and~\cite{GiannopoulouT11} and the results of Biedl and Kaufmann~\cite{BiedlK97} on optimal orthogonal  drawings of 
graphs.

The paper is organized as follows. In Section~\ref{prels} we give some basic definitions and preliminaries. In Section~\ref{prelsom}
we give a series  main combinatorial results. Based on the results of Section~\ref{prelsom}, we prove the main theorem and we derive its corollaries 
in Section~\ref{maint}. 

\newpage

\section{Preliminaries}
\label{prels}

For every positive integer $n$, let $[n]$ denote the set $\{1,2,\dots,n\}$.
A {\em graph} $G$ is a pair $(V,E)$ where $V$ is a finite set, called the {\em vertex set} 
and denoted by $V(G)$, and $E$ is a set of 2-subsets of $V$, called the {\em edge set} 
and denoted by $E(G)$. If we allow $E$ to be a multiset then $G$ is called a multigraph.
Let $G$ be a graph. For a vertex $v$, we denote by $N_G(v)$ its \emph{(open) neighborhood}, 
that is, the set of vertices which are adjacent to $v$, and by $E_{G}(v)$ the set of edges containing $v$. 
Notice that if $G$ is a multigraph $|N_{G}(v)|\leq |E_{G}(v)|$.
The degree of a vertex $v$ is $\deg_{G}(v)=|E_{G}(v)|$. 
We denote by $\Delta(G)$ the maximum degree over all vertices of $G$.

If $U\subseteq V(G)$ (respectively $u\in V(G)$ or $E\subseteq E(G)$ or $e\in E(G)$) then
$G-U$ (respectively $G-u$ or $G-E$ or $G-e$) is the graph obtained from $G$
by the removal of vertices of $U$ (respectively of vertex $u$ or edges of
$E$ or of the edge $e$). We say that a graph $H$ is a subgraph of a graph $G$, denoted by $H\subseteq G$, 
if $H$ can be obtained from $G$ after deleting edges and vertices.

We say that a graph $H$ is an {\em immersion} of a graph $G$ (or $H$ is {\em immersed} in $G$), 
$H\im G$, if there is an injective mapping $f: V(H) \to V(G)$ such that, for every edge $\{u,v\}$ of $H$, 
there is a path from $f(u)$ to $f(v)$ in $G$ and for any two distinct edges of $H$ the corresponding 
paths in $G$ are {\em edge-disjoint}, that is, they do not share common edges. 
The function $f$ is called a {\em model of $H$ in $G$}.

Let $P$ be a path and $v,u\in V(P)$. We denote by $P[v,u]$ the subpath of $P$ with endvertices $v$ and $u$.
Given two paths $P_{1}$ and $P_{2}$ who share a common endpoint $v$, we say that they are {\em well-arranged} if their common vertices appear in the same order in both paths.

A {\em tree decomposition} of a graph $G$ is a pair $(T, B)$, where $T$ is a tree and $B$ is a function 
that maps every vertex $v\in V(T)$ to a subset $B_{v}$ of $V(G)$ such that:
\begin{enumerate}[(i)]
\item $\bigcup_{v\in V(T)}B_{v}=V(G)$,

\item for every edge $e$ of $G$ there exists a vertex $t$ in $T$ such that $e \subseteq B_t$, and

\item for every $v\in V(G)$, if $r,s \in V(T)$ and $v\in B_{r} \cap B_{s}$, then for every vertex $t$ on the 
unique path between $r$ and $s$ in $T$, $v\in B_t$.
\end{enumerate}
The width of a tree decomposition $(T,B)$ is width$(T,B) := \max\{|B_{v}|-1\mid v\in V(T)\}$ and the treewidth 
of a graph $G$ is the minimum over the width$(T,B)$, where $(T,B)$ is a tree decomposition of $G$.

\paragraph{\bf Surfaces.}

A \emph{surface} $\Sigma$ is a compact 2-manifold without boundary
(we always consider connected surfaces).
Whenever we refer to a {\em
$\Sigma$-embedded graph} $G$ we consider a  2-cell embedding of
$G$ in $\Sigma$. To simplify notations, we do not distinguish
between a vertex of $G$ and the point of $\Sigma$ used in the
drawing to represent the vertex or between an edge and the line
representing it.  We also consider a graph $G$ embedded in
$\Sigma$ as the union of the points corresponding to its vertices
and edges. That way, a subgraph $H$ of $G$ can be seen as a graph
$H$, where $H\subseteq G$ in $\Sigma$.
Recall that $\Delta \subseteq \Sigma$ is
an open  (respectively closed)  disc if it is homeomorphic to
$\{(x,y):x^2 +y^2< 1\}$ (respectively $\{(x,y):x^2 +y^2\leq 1\}$).
The {\em Eüler genus} of a non-orientable surface $\Sigma$
is equal to the non-orientable genus
$\tilde{g}(\Sigma)$ (or the crosscap number).
The {\em Eüler genus}  of an orientable   surface
$\Sigma$ is $2{g}(\Sigma)$, where ${g}(\Sigma)$ is  the orientable genus
of $\Sigma$. We refer to the book of Mohar and Thomassen \cite{MoharT01} for
more details  on graphs embeddings.
The {\em Eüler genus} of a graph $G$ (denoted by $\eg(G)$) is the minimum integer $\gamma$ such 
that $G$ can be embedded on a surface of the Eüler genus $\gamma$.

\begin{figure}[h]
  \begin{center}
\scalebox{0.5}{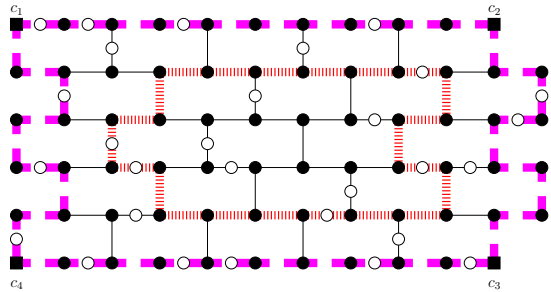}
    \end{center}
\caption{The first (magenta-dashed) and  second (red-dotted) layers of a wall of height 5}
\label{fig:layerw}
\end{figure}

\paragraph{\bf Walls.} Let $k$ and $r$ be positive integers where $k,r\geq 2$. The
\emph{$(k\times r)$-grid} $\Gamma_{k,r}$ is the Cartesian product of two paths of
lengths $k-1$ and $r-1$ respectively. 
A \emph{wall of height $k$}, $k\geq 1$, is the graph obtained from a
$((k+1)\times (2\cdot k+2))$-grid with vertices $(x,y)$,
$x\in\{1,\dots,2\cdot k+4\}$, $y\in\{1,\dots,k+1\}$, after the removal of the
``vertical'' edges $\{(x,y),(x,y+1)\}$ for odd $x+y$, and then the removal of
all vertices of degree 1. 
We denote such a wall by $W_{k}$.
The  {\em corners} of the wall $W_{k}$ are the vertices $c_{1}=(1,1)$, $c_{2}=(2\cdot k+1,1)$, $c_{3}=(2\cdot k + 1 + (k+1\mod 2),k+1)$ and $c_{4}=(1+(k+1\mod 2),k+1)$. (The square vertices in Figure~\ref{fig:layerw}.)
 
A {\em subdivided wall $W$ of height $k$} is a wall obtained from $W_{k}$ after replacing some of its 
edges by paths without common internal vertices. We call the resulting graph $W$ a {\em subdivision} 
of $W_{k}$ and the new vertices {\em subdivision vertices}. The non-subdivision vertices are called {\em original}. 
The {\em perimeter} $P$ of a subdivided wall (grid) is the cycle defined by its boundary. 

Let $W$ be a subdivided wall in a graph $G$ and $K'$ be the connected component of $G\setminus P$ 
that contains $W\setminus P$. The {\em compass} $K$ of $W$ in $G$ is the graph $G[V(K')\cup V(P)]$. 
Observe that $W$ is a subgraph of $K$ and $K$ is connected.

The {\em layers} of a subdivided wall $W$ of height $k$ are recursively defined as follows.
The first layer of $W$, denoted by $L_{1}$, is its perimeter. For $i=2,\cdots,\lceil \frac{k}{2}\rceil$, the $i$-th layer of $W$, denoted by $L_{i}$, is the 
$(i-1)$-th layer of the subwall $W'$ obtained from $W$ after removing from $W$ its perimeter and (recursively) all occurring 
vertices of degree 1 (see Figure~\ref{fig:layerw}).

We denote by $A_{i}$ the annulus defined by the cycles $L_{i}$ and $L_{i+1}$, that is, by $i$-th and $(i+1)$-th layer, $i\in [\lceil \frac{k}{2}\rceil-1]$.
Given an annulus $A$ defined by two cycles $C_{1}$ and $C_{2}$, we denote by $A^{\es}$ the interior of $A$, that is, $A\setminus (C_{1}\cup C_{2})$.

A subdivided wall of height $k$ is called {\em tight} if 
\begin{enumerate}
\item the closed disk defined by the innermost ($\lceil \frac{k}{2}\rceil$-th) layer of $W$ is edge-maximal (for reasons of uniformity we will 
denote this disk by $A_{\lceil \frac{k}{2}\rceil}$), 
\item for every $i\in [\lceil \frac{k}{2}\rceil-1]$ the annulus $A_{i}$ is edge-maximal under the condition that $A_{i+1}$ is edge-maximal.
\end{enumerate}

If $W$ is a subdivided wall of height $k$, we call {\em brick} of $W$ any facial cycle whose non-subdivided 
counterpart in $W_{h}$ has length 6. We say that two bricks are {\em neighbors} if their intersection contains an edge. 

Let $W_{k}$ be a wall. We denote by $P^{(h)}_{j}$ the shortest path connecting vertices $(1,j)$ and $(2\cdot k+2,j)$ and call these paths the {\em horizontal paths of $W_{k}$}. 
Note that these paths are vertex-disjoint.
We call the paths $P^{(h)}_{k+1}$ and 
$P^{(h)}_{1}$ the {\em southern path of $W_{k}$} and {\em northern part of $W_{k}$} respectively.

\begin{figure}[h]
\begin{center}
\scalebox{0.5}{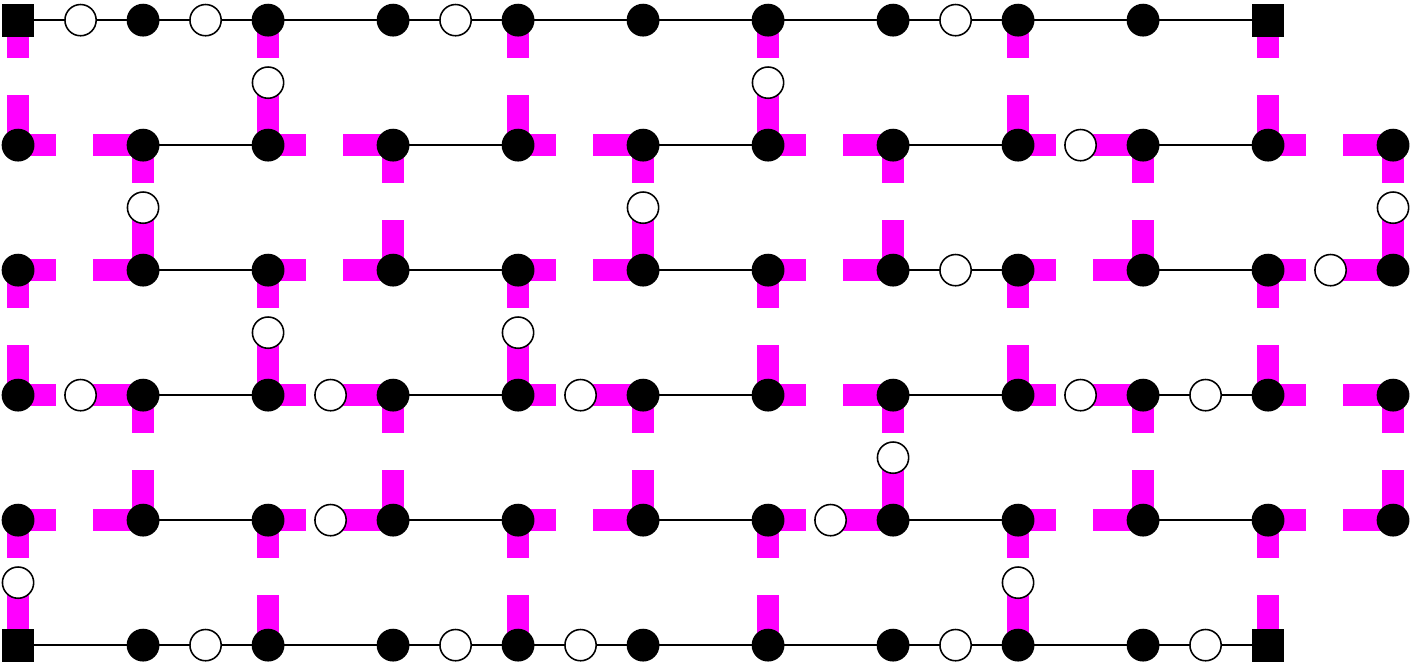}
\caption{The vertical paths of a wall of height 5}
\label{fig:vrtpths}
\end{center}
\end{figure}

Similarly, we denote by $P^{(v)}_{i}$ the shortest path connecting vertices $(i,1)$ and $(i,k+1)$ with the assumption that for, 
$i<2\cdot k+2$, $P^{(v)}_{i}$ contains only vertices $(x,y)$ with $x=i,i+1$. 
Notice that there exists a unique subfamily ${\cal P}_{v}$ of $\{P^{(v)}_{i}\mid i<2\cdot k+2\}$ of $k+1$ vertical paths with one endpoint in the southern path of $W_{k}$ and one in the northern path of $W_{k}$. 
We call these paths {\em vertical paths of $W_{k}$} and denote them by $P^{[v]}_{i}$, 
$i\in [k]$, where $P^{(v)}_{1}=P^{[v]}_{1}$ and $P^{(v)}_{2\cdot k+1}=P^{[v]}_{k+1}$. (See Figure~\ref{fig:vrtpths}.)

The paths $P^{[v]}_{1}$ and $P^{[v]}_{k+1}$ are called the {\em western part of $W_{k}$} and the 
{\em eastern part of $W_{k}$} respectively. Note that the perimeter of the wall can alternatively be 
defined as the cycle $P^{h}_{1}\cup P^{h}_{k+1}\cup P^{[v]}_{1}\cup P^{[v]}_{k+1}$.

Notice now that each vertex $u \in V(W_{k})\setminus V(P)$, is contained in exactly one vertical path, denoted by $P^{(v)}_{u}$, and in exactly one horizontal path, denoted by $P^{(h)}_{u}$, of $W_{k}$.  If $W$ is a subdivision of $W_{k}$, we will use the same notation for 
the paths obtained by the subdivisions of the corresponding paths of $W_{k}$, with further assumption that $u$ is an original vertex of $W$. \\

\begin{figure}[h]
\begin{center}
\scalebox{0.5}{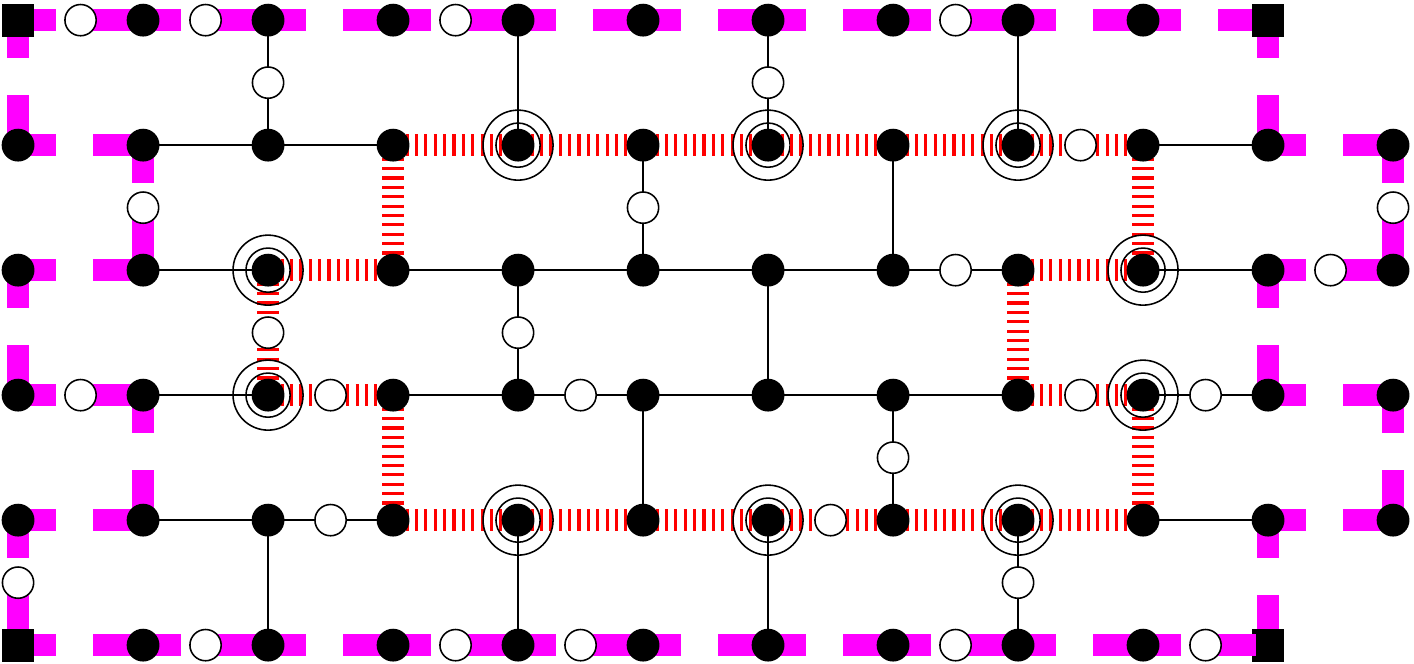}
\caption{The important vertices the second layer of a wall of height 5}
\label{fig:imptvrt}
\end{center}
\end{figure}

Given a wall $W$ and a layer $L$ of $W$, different from the perimeter of $W$. Let $W'$ be the subwall of $W$ with perimeter $L$. $W'$ is also called the {\em subwall of $W$ defined by $L$}.
We call the following vertices, {\em important} vertices of $L$; 
the original vertices of $W$ that belong to $L$ and have degree 2 in the underlying non-subdivided wall of $W'$ but are not the corners of $W'$ (where we assume that $W'$ shares the original vertices of $W$).
(See Figure~\ref{fig:imptvrt})

\begin{observation}\label{obs:impvrt}
A layer $L$ of a wall $W$ that is different from its perimeter and defines a subwall $W'$ of $W$ of height $k$ contains exactly $4k-2$ important vertices. 
\end{observation}

From Lemma 6 in~\cite{FominGT11} and Lemma 3 in~\cite{GiannopoulouT11} we obtain the following.

\begin{lemma}\label{lem:twbndwll}
Let $G$ be a graph embedded in a surface of Eüler genus $\gamma$. 
If $\tw(G)\geq 48\cdot (\gamma+1)^{\frac{3}{2}}\cdot (k+5)$, $G$ contains 
as a subgraph a subdivided wall of height $k$, whose compass in $G$ is 
embedded in a closed disk $\Delta$.
\end{lemma}

\paragraph{Confluent paths.} Let $G$ be a graph embedded in some surface $\Sigma$ 
and let $x\in V(G)$. We define a {\em disk around $x$}
as any open disk $\Delta_{x}$ with the property  that each point in $\Delta_{x}\cap G$ is  either  $x$ or belongs to the edges incident to $x$.
Let $P_{1}$ and $P_{2}$ be two edge-disjoint 
paths in $G$. 
We say that $P_{1}$ and $P_{2}$ are {\em confluent} if 
for every $x\in V(P_{1})\cap V(P_{2})$, that is not an endpoint 
of $P_{1}$ or $P_{2}$, and for every disk $\Delta_{x}$ around $x$, 
one of the connected components of 
the set $\Delta_{x}\setminus P_{1}$ does not contain any point of $P_{2}$.
We also say that a collection of paths is {\em confluent} if the paths in it are pairwise confluent. 

\begin{figure}[h]
\begin{center}
\scalebox{0.79}{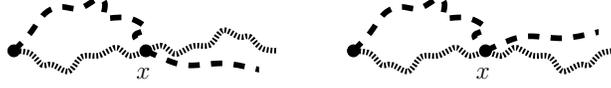}
\caption{The vertex $x$ is an overlapping vertex of the two paths on the left (dashed and dotted), while it is not an overlapping vertex of the paths on the right.}
\label{fig:cnflexmpl}
\end{center}
\end{figure}

Moreover, given two edge-disjoint paths $P_{1}$ and $P_{2}$ in $G$ we say that a vertex $x\in V(P_{1})\cap V(P_{2})$ that is not an endpoint of $P_{1}$ or $P_{2}$ is an {\em overlapping vertex of $P_{1}$ and $P_{2}$} if there exists a $\Delta_{x}$ around $x$ such that both connected components of $\Delta_{x}\setminus P_{1}$ contain points of $P_{2}$. (See, Figure~\ref{fig:cnflexmpl}.) For a family of paths ${\cal P}$, a vertex $v$ of a path $P\in {\cal P}$ is called an {\em overlapping vertex of ${\cal P}$} if there exists a path $P'\in {\cal P}$ such that $v$ is an overlapping vertex of $P$ and $P'$.

\paragraph{Orthogonal drawings.} An {\em orthogonal drawing} of a graph $G$ in a grid $\Gamma$ is a mapping which maps vertices 
$v\in V(G)$ to  subgrids $\Gamma(v)$ (called {\em boxes}) such that for every $u_{1},u_{2}\in V(G)$ 
with $u_{1}\neq u_{2}$, $\Gamma(u_{1})\cap \Gamma(u_{2})=\emptyset$, and edges 
$\{u_{1},u_{2}\}\in E(G)$ to $(u_{1}',u_{2}')$-paths whose internal vertices belong to 
$\Gamma - \bigcup_{v\in V(G)}\Gamma(v)$, their endpoints $u_{i}'$ (called {\em joining vertices of $\Gamma(u_{i})$}) belong to the perimeter of 
$\Gamma(u_{i})$, $i\in [2]$, and for every two disjoint edges $e_{i}\in E(G)$, $i\in [2]$, the corresponding 
paths are edge-disjoint.

We need the following result.
\begin{lemma}[\cite{BiedlK97}] \label{lem:grdembtamtol}
If $G$ is a simple graph then it admits an orthogonal drawing in an 
$(\frac{m+n}{2}\times \frac{m+n}{2})$-grid. Furthermore, the box size of each vertex 
$v$ is $\frac{\deg(v)+1}{2}\times \frac{\deg(v)+1}{2}$.
\end{lemma}

\section{Preliminary Combinatorial Lemmata}
\label{prelsom}

Before proving the main result of this section we first state the following lemma which we will need later on.

\begin{lemma}[\hspace{-.1mm}\cite{kurim}]
\label{tllds}
Let $r$ be a positive integer. If $G$ is a graph embedded in a surface $\Sigma$, $v,v_{1},v_{2},\dots,v_{r}\in V(G)$, and ${\cal P}$ 
is a collection of $r$ edge-disjoint paths from $v$ to $v_{1},v_{2},\dots, v_{r}$ in $G$, then $G$ contains a confluent collection ${\cal P}'$ of $r$ edge-dsjoint 
paths from $v$ to $v_{1},v_{2},\dots,v_{r}$ such that $E(\bigcup_{P\in{\cal P}'}P)\subseteq E(\bigcup_{P\in{\cal P}}P)$.
\end{lemma}

\paragraph{Detachment tree of ${\cal P}$ in $u$.}
Let $G$ be a graph embedded in a closed disk $\Delta$, $v,v_{1},v_{2},\dots,v_{k}$ be distinct vertices of $G$, and ${\cal P}=\{P_{i}\mid i\in [k]\}$ be a family of $k$ confluent edge-disjoint paths such that $P_{i}$ is a path from $v$ to $v_{i}$, $i \in[k]$. Let also $u\in V(G)\setminus \{v,v_{i}\mid i\in [k]\}$ be an internal vertex of at least two paths in ${\cal P}$. 
Let ${\cal P}_{u}= \{P_{i_{1}}$, $P_{i_{2}}$, $\dots$, $P_{i_{r}}\}$ denote the family of paths in ${\cal P}$ that contain $u$ and $\Delta_{u}$ be a disk around $u$.
Given any edge $e$ with $u\in e$ we denote by $u_{e}$ its common point with the boundary of $\Delta_{u}$.
Moreover, we denote by $e_{i_{r}}^{1}$ and $e_{i_{r}}^{2}$ the edges of $P_{i_{j}}$ incident to $u$, $j\in [r]$. 

\begin{figure}[h]
\begin{center}
\scalebox{.7}{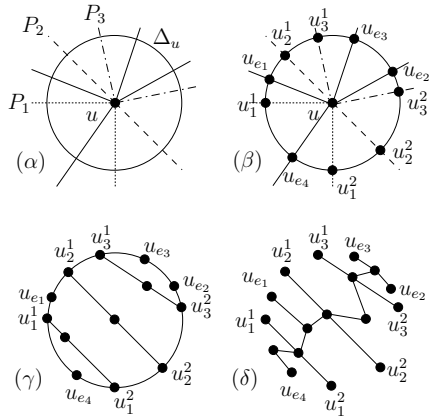}
\caption{Example of the construction of a detachment tree.}
\label{fig:dtchtr}
\end{center}

\end{figure}

We construct a tree $T_{u}$ in the following way and call it {\em detachment tree of ${\cal P}$ in $u$}.
Consider the outerplanar graph obtained from the boundary of $\Delta_{u}$ by adding the edges $\{u_{e_{i_{j}}}^{1},u_{e_{i_{j}}}^{2}\}$, $j\in [r]$.
 We subdivide the edges $\{u_{e{_{i_{j}}^{1}}},u_{e_{i_{j}}^{2}}\}$, $j\in [r]$, resulting to a planar graph. For every bounded face $f$ of the graph, let $V(f)$ denote the set of vertices that belong to $f$. We add a vertex $v_{f}$ in its interior and we make it adjacent to the vertices of $(V(f)\cap \{u_{e}\mid e\in u\})\setminus\{u_{e_{i_{j}}}^{1},u_{e_{i_{j}}}^{2}\mid j\in [r] \}$. Finally we remove the edges that lie in the boundary of $\Delta_{u}$. We call this tree $T_{u}$. Notice that for every $e$ with $u\in e$, the vertex $u_{e}$ is a leaf of $T_{u}$. (See Figure~\ref{fig:dtchtr}.)

We replace $u$ by $T_{u}$ in the following way. We first subdivide every edge $e\in G$ incident to $u$, and denote by $u_{e}$ the vertex added after the subdivision of the edge $e$. We denote by $G_{s}$ the resulting graph. 
Consider now the graph $G^{r}=(G^{s}\setminus u)\cup T_{u}$ (where, without loss of generality, we assume that $V(G\setminus u)\cap V(T_{u})=\{u_{e}\mid u\in e\}$). The graph $G^{r}$ is called {\em the graph obtained from $G$ by replacing $u$ with $T_{u}$}.
\begin{observation}\label{obs:walinvrnc}
Let $k,h$ be positive integers and $G$ be a multigraph containing as a subgraph a subdivided wall $W$ of height $h$, whose compass $C$ is embedded in a closed disk $\Delta$. Furthermore, let $v$, $v_{i}$, $i\in [k]$, be vertices of $W$ such that there exists a confluent family ${\cal P}$ of $k$ edge-disjoint paths from $v$ to the vertices $v_{i}$, $i\in [k]$. Finally, let $u\in V(C)\setminus\{v,v_{i}\mid i\in [k]\}$ belonging to more than one paths of ${\cal P}$. The graph $G^{r}$ obtained from $G$ by replacing $u$ with $T_{u}$ contains as a subgraph a subdivided wall $W'$ of height $h$, whose compass is embedded in $\Delta$ and there exists a family ${\cal P}'$ of $k$ confluent edge-disjoint paths from $v$ to $v_{i}$, $i\in [k]$, in $W'$ whose paths avoid $u$.
\end{observation}

\begin{proof}
Notice first that it is enough to prove the observation for the case where $u\in V(W)$. Let $e_{1}$, $e_{2}$ (and possibly $e_{3}$) be the edges incident to $u$ that also belong to $W$. Notice now that the vertices $u_{e_{1}}$, $u_{e_{2}}$ (and $u_{e_{3}}$) are leaves of $T_{u}$. Thus, from a folklore result, there exists a vertex $u' \in V(T_{u})$ such that there exist 2 (or 3) internally vertex-disjoint paths from $u'$ to $u_{e_{1}}$ and $u_{e_{2}}$ (and possibly $u_{e_{3}}$).
\end{proof}

We now state the following auxiliary definitions.
Let $G$ be a multigraph that contains a wall of height $k$ whose compass is embedded in a closed disk.
Let $v\in A_{\lceil \frac{k}{2}\rceil}$, that is, let $v$ be a vertex contained in the closed disk defined by the innermost layer of $W$, and let $P$ be a path from $v$ to the perimeter of $W$. For each layer $j$ of the wall, $2\leq j\leq \lceil\frac{k}{2}\rceil$, we denote by $x_{P}^{j}$ the first vertex of $P$ (starting from $v$) that also belongs to $L_{j}$ and we call it {\em incoming vertex of $P$ in $L_{j}$}.

We denote by $P^{j}$ the maximal subpath of $P$ that contains $v$ and is entirely contained in the wall defined by $L_{j}$. Moreover, we denote by $y_{P}^{j}$ its endpoint in $L_{j}$ and call it {\em outgoing vertex of $P$ in $L_{j}$}.
Notice that $x_{P}^{j}$ and $y_{P}^{j}$ are not necessarily distinct vertices.

\begin{lemma}\label{lem:vrtxord}
Let $\lambda$ and $k$ be positive integers. Let $G$ be a graph and $W$ be a tight subdivided wall of $G$ of height $k$, whose 
compass is embedded in a closed disk $\Delta$. Let also $v$ be a vertex such that $v\in A_{\lceil \frac{k}{2}\rceil}$. If there exist $\lambda$ vertex-disjoint paths 
$P_{i}$, $i\in [\lambda]$, from $v$ to vertices of the perimeter then there is a brick $B$ of $W$ with $B\cap A^{\es}_{j-1}\neq \emptyset$ that contains both $y_{P_{i}}^{j}$ and $x_{P_{i}}^{j-1}$. 
\end{lemma}

\begin{figure}[h]
\begin{center}
\scalebox{0.73}{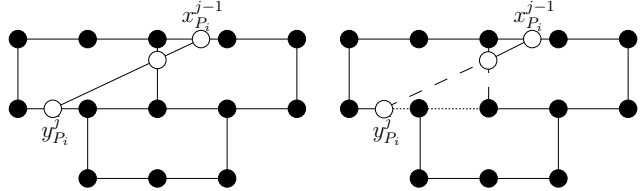}
\caption{We replace the dotted line of the wall by the dashed line.}
\label{fig:anlexmpl}
\end{center}
\end{figure}

\begin{proof}
Assume the contrary.  Then it is easy to see that we can construct an annulus $A'_{j}$ such that
$A_{j}\subsetneq A_{j}'$ and $|E(A_{j})|< |E(A_{j}')|$, a contradiction to the tightness of the wall. (See Figure~\ref{fig:anlexmpl}.)
\end{proof}

\begin{lemma}\label{vertxdispths}
Let $k$ be a positive integer and $G$ be a multigraph that contains as a subgraph a subdivided wall $W$ of height at least $4\cdot k^{2} +1$,
whose compass $K$ is embedded in a closed disk $\Delta$. Let also $V$ be a set of $k$ vertices lying in the perimeter $P$ of $W$,
whose mutual distance in the underlying non-subdivided wall is at least 2. If there exist a vertex $v\in A_{2\cdot k^{2}+1}$ and $k$ internally vertex-disjoint paths from
$v$ to vertices of $P$, then there exist $k$ internally vertex-disjoint paths from $v$ to the vertices of $V$ in $K$.
\end{lemma}

\begin{proof}
Assume first, without loss of generality, that the wall $W$ is tight. Let then $P_{1},P_{2},\dots,P_{k}$ 
be the paths from $v$ to $P$ and let $[P_{1},P_{2},\dots,P_{k},P_{1}]$ be the clockwise cyclic ordering 
according to which they appear in $W$. 
Our objective is to reroute the paths $P_{i},i\in [k]$, so that they end up to the vertices of $V$.
To do so our first step is to identify a layer of the wall for which there exist two consecutive paths whose incoming vertices on the layer are ``sufficiently far apart".

Let $j_{0}=k^{2}+1$. Consider the layer $L_{j_{0}}$ and
for every $i\in [k]$ let $T_{i}$ denote the path of $L_{j_{0}}$ starting from $x_{i}^{j_{0}}$ and ending in $x_{i+1}^{j_{0}}$ (considered clockwise), 
that is, the path of $L_{j_{0}}$ starting from the incoming vertex of $P_{i}$ in $L_{j_{0}}$ and ending to the incoming vertex of $P_{i+1}$ in $L_{j_{0}}$,
where in the case $i=k$ we abuse notation and assume that $x_{k+1}^{j_{0}}=x_{1}^{j_{0}}$ (see Figure~\ref{fig:imppths}). 
Let also $i_{0}\in [k]$ be the index such that the path $T_{i_{0}}$ 
contains the maximum number of important vertices amongst the $T_{i}$'s. 
Without loss of generality we may assume that $i_{0}=\lceil\frac{k}{2}\rceil$. From Observation~\ref{obs:impvrt}, as $L_{j_{0}}$ defines a subwall of $W$ of height $2 k^{2}+1$, $L_{j_{0}}$ contains exactly $8k^{2}+2$ important vertices.
Thus, at least $7k$ important vertices are internally contained in $T_{i_{0}}$.
This concludes the first step of the proof.

\begin{figure}[h]
\begin{center}
\scalebox{0.8}{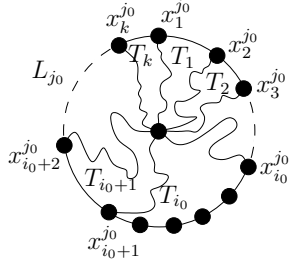}
\caption{The $T_{i}$'s, $i\in [k]$}
\label{fig:imppths}
\end{center}
\end{figure}

Let now $j_{1}=k+1$. At the next step, using the part of the wall that is contained in $A[L_{j_{0}},L_{j_{1}}]$, that is, 
in the annulus between the $j_{0}$-th and the $j_{1}$-th layer of the wall, 
we find $k$ internally vertex-disjoint paths from the incoming vertices of the paths in $L_{j_{0}}$ to $k$ consecutive important vertices of the
$k+1$-th layer of the wall. These are the paths that will allow us to reroute the original paths.

 Continuing the proof, let $u_{1},u_{2},\dots,u_{k}$ be a set of successive important vertices appearing clockwise in $T_{i_{0}}$ 
 such that the paths $T_{i_{0}}[x_{i_{0}}^{j_{0}},u_{1}]$ and $T_{i_{0}}[u_{k},x_{i_{0}+1}^{j_{0}}]$ internally contain at least $3k$ important vertices.
Notice that, without loss of generality, we may assume that the vertices $u_{i}$, $i\in [k]$, belong to the northern part of $W'$.
Recall here that each original vertex $w$ of $W'\subseteq W\setminus P$ is contained in exactly one vertical path $P^{(v)}_{w}$ of $W$.
For every $i\in [k]$ we assign the path $R_{i}$ to the vertex $u_{i}$ in the following way. Let $R_{i}$ be the maximal subpath of $P_{u_{i}}^{(v)}$ whose endpoints are $u_{i}$ and the important vertex of $L_{j_{1}}$ that also belongs to $P_{u_{i}}^{(v)}$, which from now on we will denote by $u_{i}^{f}$.
Note here that the paths $R_{i}$, $i\in [k]$, are vertex-disjoint 
and  do not contain any of the vertices belonging to the interior of the disk defined by $L_{j_{0}}$ in the compass of $W$ 
(See, for example, Figure~\ref{fig:pthsexmpl}).

\begin{figure}[h]
\begin{center}
\scalebox{0.45}{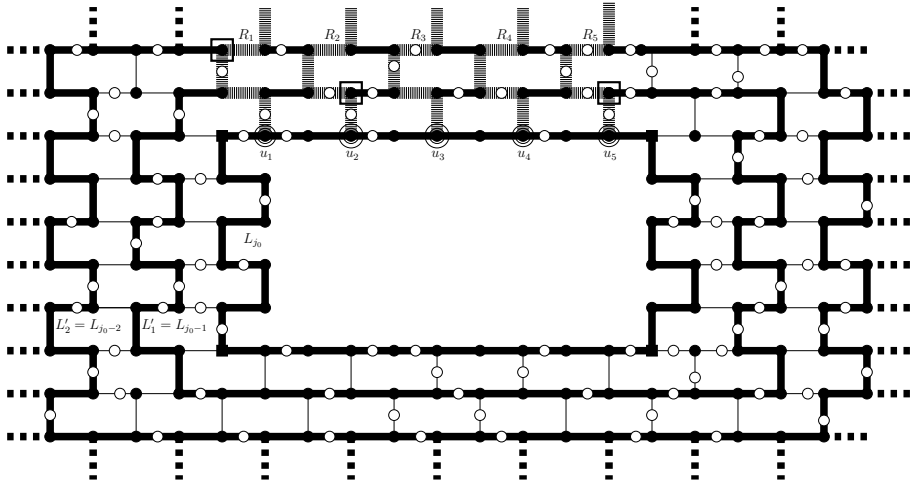}
\caption{The important vertices of $L_{j_{0}}$, the layers $L_{1}'$ and $L_{2}'$, and the paths $R_{i}$}
\label{fig:pthsexmpl}
\end{center}
\end{figure}

\noindent Notice now that $T_{i_{0}}$, and thus $L_{j_{0}}$, contains a path $F_{1}$ from $x_{i_{0}}^{j_{0}}$ to $u_{i_{0}}$ and a path $F_{2}$ from $u_{i_{0}+1}$ 
to $x_{i_{0}+1}^{j_{0}}$ that are vertex-disjoint and do not contain vertices of any path other than $P_{i_{0}}$ and $P_{i_{0}+1}$.
Consider now the $\lceil \frac{k-2}{2}\rceil$ consecutive layers of $W$ preceeding $L_{j_{0}}$, that is, the layers $L_{j}'=L_{j_{0}-j}$, $j\in [\lceil \frac{k-2}{2}\rceil]$. For every $j\in [\lceil \frac{k-2}{2}\rceil]$ let $u_{i_{0}-j}^{j}$
be the first time the path $R_{i_{0} -j}$ meets $L_{j}'$ starting from $u_{i_{0} -j}$
and $u_{i_{0}+1+j}^{j}$ be the first time the path $R_{i_{0}+1+j}$
meets $L_{j}'$ starting from $u_{i_{0}+1+j}$. (See, for example, the vertices inside the squares in Figure~\ref{fig:pthsexmpl}.)

We need to prove the following. \\

\noindent{\em Claim:}
For every $j\in [\lceil\frac{k-2}{2}\rceil]$, there exist two vertex-disjoint paths $F^{1}_{j}$ and $F^{2}_{j}$ between the pairs of vertices $(x_{i_{0}-j}^{j_{0}-j},u_{i_{0}-j}^{j})$ and $(x_{i_{0}+1+j}^{j_{0}-j},u_{i_{0}+1+j}^{j})$ that do not intersect the paths $\{R_{l}\mid i_{0}-j < l< i_{0}+1+j\}$.\\

\noindent{\em Proof of Claim:}
Indeed, this holds by inductively applying the combination of Lemma~\ref{lem:vrtxord} with the assertion that
for every $j\leq 2\cdot k^{2}+1$ and every $p,q$ with $1<p<q<k$, the outgoing vertices of $P_{p-1}$ and $P_{q+1}$ and the incoming vertices of $P_{p}$ and $P_{q}$ in the layer $L_{j}$, $y_{p-1}^{j}$, $y_{q+1}^{j}$, $x_{p}^{j}$, and $x_{q}^{j}$ respectively appear in $L_{j}$ respecting the clockwise order
$$[y_{p-1}^{j},x_{p}^{j},x_{q}^{j},y_{q+1}^{j}]$$
in the tight wall $W$. This completes the proof of the claim.\hfill$\diamond$\\

\noindent We now construct the following paths.
First, let $$Q_{i_{0}}=F_{1}\cup R_{i_{0}}$$ and $$Q_{i_{0}+1}=F_{2}\cup R_{i_{0}+1},$$ that is, $Q_{i_{0}}$ is the union of the paths $F_{1}$ and $R_{i_{0}}$, and $Q_{i_{0}+1}$ is the union of the paths $F_{2}$ and $R_{i_{0}+1}$.
Then, for every $j\in [\lceil \frac{k-2}{2}\rceil]$, let $$Q_{i_{0}-j}=P_{i_{0}-j}[x_{i_{0}-j}^{j_{0}},x_{i_{0}-j}^{j_{0}-j}]\cup F_{j}^{1}\cup R_{i_{0} -j}[u_{i_{0}-j}^{j},u_{i_{0}-j}^{f}],$$ that is, $Q_{j_{0}-j}$ is the union of the following paths; (a) the subpath of $P_{i_{0}-j}$ between its incoming vertex in the $j_{0}$-th layer and its incoming vertex in the $(j_{0}-j)$-th layer, (b) the path $F_{j}^{1}$ defined in the claim above, and (c) the subpath of $R_{i_{0}-j}$ between the vertices $u_{i_{0}-j}^{j}$ and $u_{i_{0}-j}^{f}$.

Finally, for every $j\in [\lceil \frac{k-2}{2}\rceil]$ ($j\in [\lceil \frac{k-2}{2}\rceil-1]$, if $k$ is odd) let $$Q_{i_{0}+1+j}=P_{i_{0}+1+j}[x_{i_{0}+1+j}^{j_{0}},x_{i_{0}+1+j}^{j_{0}-j}]\cup F_{j}^{2}\cup R_{i_{0} +1+j}[u_{i_{0}+1+j}^{j},u_{i_{0}+1+j}^{f}],$$
that is, $Q_{i_{0}+1+j}$ is the union of the following three paths; (a) the subpath of $P_{i_{0}+1+j}$ between its incoming vertex in the $j_{0}$-th layer and its incoming vertex in the $(j_{0}-j)$-th layer, (b) the path $F_{j}^{2}$ defined in the  claim above, and (c) the subpath of $R_{i_{0}+1+j}$ between the vertices $u_{i_{0}+1+j}^{j}$ and $u_{i_{0}+1+j}^{f}$.

From the claim above and Lemma~\ref{lem:vrtxord} we get that the above paths are vertex-disjoint. This concludes the second step of the proof.

We claim now that we may reroute the paths $P_{i}$, $i\in [k]$, in such a way that they end up to the vertices $u_{i}^{f}$, $i\in [k]$. Indeed, let $P_{i}'=P_{i}[v,x_{i}^{j_{0}}]\cup Q_{i}$, $i\in [k]$. From their construction these paths are vertex-disjoint and end up to the vertices $u_{i}^{f}$, $i\in [k]$. (For a rough estimation of the position of the paths in the wall see Figure~\ref{fig:pthsexmpl2}.)

\begin{figure}[h]
\begin{center}
\scalebox{0.45}{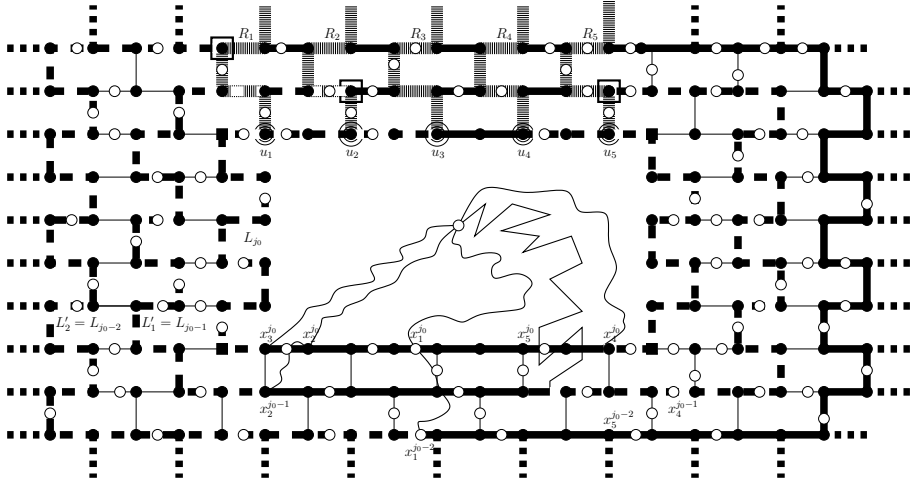}
\caption{Part of the rerouted paths}
\label{fig:pthsexmpl2}
\end{center}
\end{figure}

Concluding the proof, 
as the mutual distance of the vertices of $V$ in the underlying non-subdivided wall is at least 2, it is easy to notice that in the annulus defined by $L_{j_{1}}$ and $L_{1}$ there exist $k$ vertex-disjoint paths from the vertices $u_{i}^{f}$, $i\in [k]$, to the vertices of $V$.
\end{proof}

We now prove the main result of this section.

\begin{lemma}\label{lem:fixhighrdeg}
Let $k$ be a positive integer and $G$ be a $k$-edge-connected multigraph embedded in a surface of 
Eüler genus $\gamma$ that contains a subdivided wall $W$ of height at least $4\cdot k^{2}+1$ as a subgraph, whose 
compass $C$ is embedded in a closed disk $\Delta$. Let also $S$ be a set of vertices in the 
perimeter of $W$ whose mutual distance in the underlying non-subdivided wall is at least 2. If $|S|\leq k$ then 
there exist a vertex $v$ in $W$ and $|S|$ edge-disjoint paths from $v$ to the vertices of $S$.
\end{lemma}

\begin{proof}
Let $v\in A_{2k^{2}+1}$ and $u\in L_{1}$ be vertices belonging to the closed disk defined by the layer $L_{2 k^{2}+1}$ and to the perimeter of the wall respectively.
As $G$ is $k$-edge-connected there exist $k$ edge-disjoint paths $P_{1},P_{2},\dots,P_{k}$ connecting $v$ and $u$.
By Lemma~\ref{tllds}, we may assume that the paths are confluent. Let ${\cal P}'=\{P_{i}'\mid i\in [k]\}$ be the family of paths $P_{i}'=P_{i}[v,x_{i}^{1}]$, $i\in [k]$,
that is, let ${\cal P}'$ be the family of paths consisting of he subpaths of $P_{i}$, $i\in [k]$, between $v$ and the first vertex on which they meet the perimeter of $W$.

Let $V$ be the set of vertices in $V(C)\setminus (V(L_{1})\cup\{v\})$ that are contained in more than one path in ${\cal P}'$. We obtain the graph $\hat{G}$ by replacing every vertex $z\in V$ with the detachment tree
of ${\cal P}'$ in $z$. From Observation~\ref{obs:walinvrnc}, $\hat{G}$ contains a wall $\hat{W}$ of height $4k^{2}+1$ whose compass is embedded in $\Delta$. Notice also that, as no changes have occurred in the perimeter of $W$, $W$ and $\hat{W}$ share the same perimeter. Furthermore, $\hat{W}$ contains $k$ internally vertex-disjoint paths from $v$ to the perimeter of $\hat{W}$.  Thus, from Lemma~\ref{vertxdispths}, $\hat{W}$ contains $k$ vertex-disjoint paths from $v$ to $S$. It is now easy to see, by contracting each one of the trees $T_{z}$, $z\in V(C)\setminus (V(L_{1})\cup\{v\})$, to a single vertex that $W$ contains $k$ edge-disjoint paths from $v$ to $S$.
\end{proof}

\section{Main Theorem}
\label{maint}

By combining Lemmata~\ref{lem:fixhighrdeg},~\ref{lem:twbndwll} and~\ref{lem:grdembtamtol} we obtain the following.

\begin{theorem}\label{mainthm}
There exists a computable function $f:\mathbb{N}\rightarrow \mathbb{N}$ such that for every 
multigraph $G$ of Eüler genus $\gamma$ and every connected graph $H$ one of the following holds:
\begin{enumerate}
\item $\tw(G)\leq f(\gamma)\cdot \lambda \cdot k$, where $\lambda=\Delta(H)$ and $k=\m(H)$
\item $G$ is not $\lambda$-edge-connected,
\item $H\im G$.
\end{enumerate}
\end{theorem}

\begin{proof}
Let 
$$f(\gamma,\lambda,k)=48\cdot \left(\gamma+1\right)^{\frac{3}{2}}\cdot \left(\frac{4\left(4\lambda+1\right)k}{2}+5\right),$$ 

\noindent and assume that $\tw(G)\geq f(\gamma,\lambda,k)$ and $G$ is $\lambda$-edge-connected.
From Lemma~\ref{lem:twbndwll}, we obtain that $G$ contains as a subgraph a subdivided wall $W$ of height $2 (2\lambda +1)k$ whose compass is embedded in a closed disk.

In what follows we will construct a model of $H$ into the wall.
From Lemma~\ref{lem:grdembtamtol}, $H$ admits an a orthogonal drawing $\psi$ in an $$\left(\frac{\m(H)+\n(H)}{2}\times \frac{\m(H)+\n(H)}{2}\right)\text{-grid},$$ where the box size of each vertex $v\in V(H)$ is $\frac{\deg(v)+1}{2}\times\frac{\deg(v)+1}{2}$.

 Notice now that $\psi$ can be scaled to an orthogonal drawing $\phi$ to the grid $\Gamma$ of size $$\left(\frac{2\left(4\lambda+1\right)\left(\m(H)+\n(H)\right)}{2}+1\right)\times 2\left(\frac{2\left(4\lambda+1\right)\left(\m(H)+\n(H)\right)+2}{2}+1\right),$$
where the box size of each vertex is $(4(\deg(v))^{2}+2)\times 2(4(\deg(v))^{2}+2)$, the joining vertices of each box have mutual distance at least 2 in the perimeter of the box and no joining vertex is a corner of the box.

Moreover, for every vertex $u$, $u\in \img(\phi)\setminus \cup_{v\in V(H)}\Gamma(v)$ of degree $4$, that is, for every vertex in the image of $\phi$ that is the intersection of two paths, there is a box in the grid of size $(4\deg(u)^{2}+2)\times 2(4\deg(u)^{2}+2)$, denoted by $Q(u)$, containing only this vertex and vertices of the paths it belongs to. 
We denote by $u^{i}$, $i\in [4]$, the vertices of $\img(\phi)$ belonging to the boundary of $Q(u)$ and, for uniformity, also call them {\em joining vertices of $Q(u)$}.

Towards finding a model of $H$ in the wall observe
that the grid $\Gamma$ contains as a subgraph a wall of height $\left(4\lambda+1\right)\left(\m(H)+\n(H)\right)$
such that each one of the boxes, either $\Gamma(v)$, $v\in V(H)$, or $Q(v)$, where $v$ is the intersection of two paths in the image of $\phi$ contains a wall $W(v)$ of height $4\deg(v)^{2}+1$ and the joining vertices of $\Gamma(v)$ (the vertices $v^{i}$, $i\in [4]$, respectively) belong to the perimeter of the wall and have distance at least 2 in it. Consider now the mapping of $H$ to $W$ where the boxes $\Gamma(v)$ and $Q(v)$ are mapped into subwalls $W(v)$ of $W$ of height $4\deg(v)^{2}+1$ joined together by vertex-disjoint paths as given by the orthogonal drawing $\phi$.
 From Lemma~\ref{lem:fixhighrdeg}, as every $W(v)$ has height $4\deg(v)^{2}+1$ and its compass is embedded in a closed disk, there exist a vertex $z_{v}\in V(W(v))$ and $\deg(v)$ edge-disjoint paths from $z_{v}$ to the joining vertices of $W(v)$. 
 It is now easy to see that $W$ contains a model of $H$.
\end{proof}

\noindent Notice now that in the case when $\Delta(H)=O(1)$ we get the following.

\begin{theorem}\label{bdthm}
There exists a computable function $f:\mathbb{N}\rightarrow \mathbb{N}$ such that for every 
multigraph $G$ of Eüler genus $\gamma$ and every connected graph $H$ one of the following holds:
\begin{enumerate}
\item $\tw(G)\leq f(\gamma) \cdot  \n(H)$, 
\item $G$ is not $\Delta(H)$-edge-connected,
\item $H\im G$.
\end{enumerate}
\end{theorem}

\noindent The following two corollaries are immediate consequences of Theorems~\ref{mainthm} and~\ref{bdthm}.

\begin{corollary}
There exists a computable function $f:\mathbb{N}\rightarrow \mathbb{N}$ such that 
for every multigraph $G$ of Eüler genus $\gamma$ and every $k\in \mathbb{N}$ one of the following holds:
\begin{enumerate}
\item $\tw(G)\leq f(\gamma)\cdot k^{3}$,
\item $G$ is not $k$-edge-connected,
\item $K_{k+1} \im G$.
\end{enumerate}
\end{corollary}

\begin{corollary}
There exists a computable function $f:\mathbb{N}\rightarrow \mathbb{N}$ such that 
for every multigraph $G$ of Eüler genus $\gamma$ and every $k\in \mathbb{N}$ one of the following holds:
\begin{enumerate}
\item $\tw(G)\leq f(\gamma)\cdot k^{2}$,
\item $G$ is not $4$-edge-connected,
\item $(k\times k)$-grid is an immersion of $G$.
\end{enumerate}
\end{corollary}

\noindent However, when $H$ is the grid a straightforward argument gives the following result.

\begin{theorem}\label{thm:exclgrdimrs}
There exists a computable function $f:\mathbb{N}\rightarrow \mathbb{N}$ such that 
for every multigraph $G$ that is embedded in a surface of Eüler genus $\gamma$ and every $k\in \mathbb{N}$ one of the following holds:
\begin{enumerate}
\item $\tw(G)\leq f(\gamma)\cdot k$.
\item $G$ is not $4$-edge-connected.
\item $(k\times k)$-grid is an immersion of $G$.
\end{enumerate}
\end{theorem}

\begin{proof}
Let
 $$f(\gamma,k)=48\cdot (\gamma+1)^{\frac{3}{2}}\cdot ((4^{3}+3)\cdot k+5).$$

\noindent Assume that $G$ is $4$-edge-connected and that $\tw(G)\geq f(\gamma,k)$.
As $\tw(G)\geq f(\gamma,k)$, from Lemma~\ref{lem:twbndwll} it follows that $G$ contains as a subgraph
a subdivided wall $W$ of height $(4^{3}+3)k$, whose compass in $G$ is embedded in a closed disk $\Delta$.

Consider the $k^{2}$ subwalls of $W$ of height $(4^{3}+1)$ that occur after removing from it the paths
$P_{(4^{3}+3)j}^{[v]}$, $P_{(4^{3}+3)j}^{[h]}$, $i,j\in [k]$.
For every $i,j\in [k]$, we denote by $W_{(i,j)}$ the subwall that is contained inside the disk that is defined by the paths
$P^{(h)}_{(4^{3}+3)(i-1)}$, $P^{(h)}_{(4^{3}+3)i}$, $P^{[v]}_{(4^{3}+3)(j-1)}$, and $P^{[v]}_{(4^{3}+3)j}$.
In the case where $j=1$ and $i=1$, we abuse notation and consider as $P^{(h)}_{(4^{3}+3)(j-1)}$ and $P^{[v]}_{(4^{3}+3)(j-1)}$ the paths $P^{(h)}_{1}$ and $P^{[v]}_{1}$, respectively.

 From Lemma~\ref{lem:fixhighrdeg} and the hypothesis that $G$ is $4$-edge-connected, for $k=4$, it follows that
in the compass of each one of the subwalls $\{W_{(i,j)}\mid i,j\in [k]\}$ we may find a vertex $v_{(i,j)}$ and four edge-disjoint paths from $v_{(i,j)}$ 
to the vertices $v_{(i,j)}^{n}$, $v_{(i,j)}^{s}$, $v_{(i,j)}^{w}$, and $v_{(i,j)}^{e}$, that lie in the
northern, southern, western, and eastern path of the wall, respectively.

Finally, we consider the function $g((i,j))=v_{(i,j)}$ that maps the vertex $(i,j)$ of the $(k\times k)$-grid to the vertex $v_{(i,j)}$ of the wall $W_{(i,j)}$.
Is now easy to see that $g$ is an immersion model of the $(k\times k)$-grid in the compass of the wall $W$ and the theorem follows as $f$ is linear on $k$.
\end{proof}

\section{Conclusions}

In this paper, we proved sufficient conditions for  the containment of any connected graph $H$ as an immersion in graphs of bounded genus.
We would like to remark here that our proofs also hold if we, instead, consider the strong immersion 
relation where we additionally ask that the paths of the model $f$ of $H$ in $G$ that 
correspond to the edges of $H$ are internally disjoint from $f(V(H))$.

In our results, it appears that both big treewidth and 
the edge connectivity requirement are necessary in order to enforce the appearance of
a graph as an immersion. A natural open problem to investigate is the existence 
of counterparts of  our results for the case of the topological minor relation.
Certainly, here edge connectivity should be replaced by vertex connectivity. However, what we can only report 
 is that stronger conditions than 
just asking for sufficiently big treewidth are required for such an extension.

\bibliography{complete}
\bibliographystyle{abbrv}

\end{document}

%% file: 1.pdf_t
\begin{picture}(0,0)%
\includegraphics{1.pdf}%
\end{picture}%
\setlength{\unitlength}{3947sp}%
\begingroup\makeatletter\ifx\SetFigFont\undefined%
\gdef\SetFigFont#1#2#3#4#5{%
  \reset@font\fontsize{#1}{#2pt}%
  \fontfamily{#3}\fontseries{#4}\fontshape{#5}%
  \selectfont}%
\fi\endgroup%
\begin{picture}(6776,3711)(361,-2989)
\put(376,539){\makebox(0,0)[lb]{\smash{{\SetFigFont{12}{14.4}{\rmdefault}{\mddefault}{\updefault}{\color[rgb]{0,0,0}$c_{1}$}%
}}}}
\put(6376,539){\makebox(0,0)[lb]{\smash{{\SetFigFont{12}{14.4}{\rmdefault}{\mddefault}{\updefault}{\color[rgb]{0,0,0}$c_{2}$}%
}}}}
\put(6376,-2911){\makebox(0,0)[lb]{\smash{{\SetFigFont{12}{14.4}{\rmdefault}{\mddefault}{\updefault}{\color[rgb]{0,0,0}$c_{3}$}%
}}}}
\put(376,-2911){\makebox(0,0)[lb]{\smash{{\SetFigFont{12}{14.4}{\rmdefault}{\mddefault}{\updefault}{\color[rgb]{0,0,0}$c_{4}$}%
}}}}
\end{picture}%

%% file: vrtpths.pdf_t
\begin{picture}(0,0)%
\includegraphics{vrtpths.pdf}%
\end{picture}%
\setlength{\unitlength}{3947sp}%
\begingroup\makeatletter\ifx\SetFigFont\undefined%
\gdef\SetFigFont#1#2#3#4#5{%
  \reset@font\fontsize{#1}{#2pt}%
  \fontfamily{#3}\fontseries{#4}\fontshape{#5}%
  \selectfont}%
\fi\endgroup%
\begin{picture}(6773,3174)(364,-2698)
\end{picture}%

%% file: impvrt.pdf_t
\begin{picture}(0,0)%
\includegraphics{impvrt.pdf}%
\end{picture}%
\setlength{\unitlength}{3947sp}%
\begingroup\makeatletter\ifx\SetFigFont\undefined%
\gdef\SetFigFont#1#2#3#4#5{%
  \reset@font\fontsize{#1}{#2pt}%
  \fontfamily{#3}\fontseries{#4}\fontshape{#5}%
  \selectfont}%
\fi\endgroup%
\begin{picture}(6773,3174)(364,-2698)
\end{picture}%

%% file: cnflexmpl.pdf_t
\begin{picture}(0,0)%
\includegraphics{cnflexmpl.pdf}%
\end{picture}%
\setlength{\unitlength}{3947sp}%
\begingroup\makeatletter\ifx\SetFigFont\undefined%
\gdef\SetFigFont#1#2#3#4#5{%
  \reset@font\fontsize{#1}{#2pt}%
  \fontfamily{#3}\fontseries{#4}\fontshape{#5}%
  \selectfont}%
\fi\endgroup%
\begin{picture}(4899,754)(396,-55)
\put(4126, 14){\makebox(0,0)[lb]{\smash{{\SetFigFont{12}{14.4}{\rmdefault}{\mddefault}{\updefault}{\color[rgb]{0,0,0}$x$}%
}}}}
\put(1426, 14){\makebox(0,0)[lb]{\smash{{\SetFigFont{12}{14.4}{\rmdefault}{\mddefault}{\updefault}{\color[rgb]{0,0,0}$x$}%
}}}}
\end{picture}%

%% file: dtchtr2.pdf_t
\begin{picture}(0,0)%
\includegraphics{dtchtr2.pdf}%
\end{picture}%
\setlength{\unitlength}{3947sp}%
\begingroup\makeatletter\ifx\SetFigFont\undefined%
\gdef\SetFigFont#1#2#3#4#5{%
  \reset@font\fontsize{#1}{#2pt}%
  \fontfamily{#3}\fontseries{#4}\fontshape{#5}%
  \selectfont}%
\fi\endgroup%
\begin{picture}(3695,3763)(73,-2748)
\put(3041,-729){\makebox(0,0)[lb]{\smash{{\SetFigFont{12}{14.4}{\rmdefault}{\mddefault}{\updefault}{\color[rgb]{0,0,0}$u^{2}_{1}$}%
}}}}
\put( 88, 30){\makebox(0,0)[lb]{\smash{{\SetFigFont{12}{14.4}{\rmdefault}{\mddefault}{\updefault}{\color[rgb]{0,0,0}$P_{1}$}%
}}}}
\put(789,-70){\makebox(0,0)[lb]{\smash{{\SetFigFont{12}{14.4}{\rmdefault}{\mddefault}{\updefault}{\color[rgb]{0,0,0}$u$}%
}}}}
\put(1390,613){\makebox(0,0)[lb]{\smash{{\SetFigFont{12}{14.4}{\rmdefault}{\mddefault}{\updefault}{\color[rgb]{0,0,0}$\Delta_{u}$}%
}}}}
\put(213,724){\makebox(0,0)[lb]{\smash{{\SetFigFont{12}{14.4}{\rmdefault}{\mddefault}{\updefault}{\color[rgb]{0,0,0}$P_{2}$}%
}}}}
\put(767,832){\makebox(0,0)[lb]{\smash{{\SetFigFont{12}{14.4}{\rmdefault}{\mddefault}{\updefault}{\color[rgb]{0,0,0}$P_{3}$}%
}}}}
\put(3539,-411){\makebox(0,0)[lb]{\smash{{\SetFigFont{12}{14.4}{\rmdefault}{\mddefault}{\updefault}{\color[rgb]{0,0,0}$u^{2}_{2}$}%
}}}}
\put(2189,435){\makebox(0,0)[lb]{\smash{{\SetFigFont{12}{14.4}{\rmdefault}{\mddefault}{\updefault}{\color[rgb]{0,0,0}$u_{e_{1}}$}%
}}}}
\put(3259,704){\makebox(0,0)[lb]{\smash{{\SetFigFont{12}{14.4}{\rmdefault}{\mddefault}{\updefault}{\color[rgb]{0,0,0}$u_{e_{3}}$}%
}}}}
\put(2778,-62){\makebox(0,0)[lb]{\smash{{\SetFigFont{12}{14.4}{\rmdefault}{\mddefault}{\updefault}{\color[rgb]{0,0,0}$u$}%
}}}}
\put(1308,-1247){\makebox(0,0)[lb]{\smash{{\SetFigFont{12}{14.4}{\rmdefault}{\mddefault}{\updefault}{\color[rgb]{0,0,0}$u_{e_{3}}$}%
}}}}
\put(3636,333){\makebox(0,0)[lb]{\smash{{\SetFigFont{12}{14.4}{\rmdefault}{\mddefault}{\updefault}{\color[rgb]{0,0,0}$u_{e_{2}}$}%
}}}}
\put(2147, 13){\makebox(0,0)[lb]{\smash{{\SetFigFont{12}{14.4}{\rmdefault}{\mddefault}{\updefault}{\color[rgb]{0,0,0}$u^{1}_{1}$}%
}}}}
\put(158,-496){\makebox(0,0)[lb]{\smash{{\SetFigFont{12}{14.4}{\rmdefault}{\mddefault}{\updefault}{\color[rgb]{0,0,0}$(\alpha)$}%
}}}}
\put(2104,-488){\makebox(0,0)[lb]{\smash{{\SetFigFont{12}{14.4}{\rmdefault}{\mddefault}{\updefault}{\color[rgb]{0,0,0}$(\beta)$}%
}}}}
\put(143,-2449){\makebox(0,0)[lb]{\smash{{\SetFigFont{12}{14.4}{\rmdefault}{\mddefault}{\updefault}{\color[rgb]{0,0,0}$(\gamma)$}%
}}}}
\put(2093,-2445){\makebox(0,0)[lb]{\smash{{\SetFigFont{12}{14.4}{\rmdefault}{\mddefault}{\updefault}{\color[rgb]{0,0,0}$(\delta)$}%
}}}}
\put(2554,-620){\makebox(0,0)[lb]{\smash{{\SetFigFont{12}{14.4}{\rmdefault}{\mddefault}{\updefault}{\color[rgb]{0,0,0}$u_{e_{4}}$}%
}}}}
\put(602,-2571){\makebox(0,0)[lb]{\smash{{\SetFigFont{12}{14.4}{\rmdefault}{\mddefault}{\updefault}{\color[rgb]{0,0,0}$u_{e_{4}}$}%
}}}}
\put(1558,-2397){\makebox(0,0)[lb]{\smash{{\SetFigFont{12}{14.4}{\rmdefault}{\mddefault}{\updefault}{\color[rgb]{0,0,0}$u^{2}_{2}$}%
}}}}
\put(177,-1900){\makebox(0,0)[lb]{\smash{{\SetFigFont{12}{14.4}{\rmdefault}{\mddefault}{\updefault}{\color[rgb]{0,0,0}$u^{1}_{1}$}%
}}}}
\put(183,-1671){\makebox(0,0)[lb]{\smash{{\SetFigFont{12}{14.4}{\rmdefault}{\mddefault}{\updefault}{\color[rgb]{0,0,0}$u_{e_{1}}$}%
}}}}
\put(1659,-1560){\makebox(0,0)[lb]{\smash{{\SetFigFont{12}{14.4}{\rmdefault}{\mddefault}{\updefault}{\color[rgb]{0,0,0}$u_{e_{2}}$}%
}}}}
\put(1040,-2655){\makebox(0,0)[lb]{\smash{{\SetFigFont{12}{14.4}{\rmdefault}{\mddefault}{\updefault}{\color[rgb]{0,0,0}$u^{2}_{1}$}%
}}}}
\put(834,-1170){\makebox(0,0)[lb]{\smash{{\SetFigFont{12}{14.4}{\rmdefault}{\mddefault}{\updefault}{\color[rgb]{0,0,0}$u^{1}_{3}$}%
}}}}
\put(1732,-1791){\makebox(0,0)[lb]{\smash{{\SetFigFont{12}{14.4}{\rmdefault}{\mddefault}{\updefault}{\color[rgb]{0,0,0}$u^{2}_{3}$}%
}}}}
\put(2129,-1927){\makebox(0,0)[lb]{\smash{{\SetFigFont{12}{14.4}{\rmdefault}{\mddefault}{\updefault}{\color[rgb]{0,0,0}$u^{1}_{1}$}%
}}}}
\put(2786,-1170){\makebox(0,0)[lb]{\smash{{\SetFigFont{12}{14.4}{\rmdefault}{\mddefault}{\updefault}{\color[rgb]{0,0,0}$u^{1}_{3}$}%
}}}}
\put(504,-1293){\makebox(0,0)[lb]{\smash{{\SetFigFont{12}{14.4}{\rmdefault}{\mddefault}{\updefault}{\color[rgb]{0,0,0}$u^{1}_{2}$}%
}}}}
\put(2189,-1516){\makebox(0,0)[lb]{\smash{{\SetFigFont{12}{14.4}{\rmdefault}{\mddefault}{\updefault}{\color[rgb]{0,0,0}$u_{e_{1}}$}%
}}}}
\put(2455,-1293){\makebox(0,0)[lb]{\smash{{\SetFigFont{12}{14.4}{\rmdefault}{\mddefault}{\updefault}{\color[rgb]{0,0,0}$u^{1}_{2}$}%
}}}}
\put(2529,-2541){\makebox(0,0)[lb]{\smash{{\SetFigFont{12}{14.4}{\rmdefault}{\mddefault}{\updefault}{\color[rgb]{0,0,0}$u_{e_{4}}$}%
}}}}
\put(2931,-2670){\makebox(0,0)[lb]{\smash{{\SetFigFont{12}{14.4}{\rmdefault}{\mddefault}{\updefault}{\color[rgb]{0,0,0}$u^{2}_{1}$}%
}}}}
\put(3489,-2367){\makebox(0,0)[lb]{\smash{{\SetFigFont{12}{14.4}{\rmdefault}{\mddefault}{\updefault}{\color[rgb]{0,0,0}$u^{2}_{2}$}%
}}}}
\put(3494,-1959){\makebox(0,0)[lb]{\smash{{\SetFigFont{12}{14.4}{\rmdefault}{\mddefault}{\updefault}{\color[rgb]{0,0,0}$u^{2}_{3}$}%
}}}}
\put(3599,-1636){\makebox(0,0)[lb]{\smash{{\SetFigFont{12}{14.4}{\rmdefault}{\mddefault}{\updefault}{\color[rgb]{0,0,0}$u_{e_{2}}$}%
}}}}
\put(3109,-1202){\makebox(0,0)[lb]{\smash{{\SetFigFont{12}{14.4}{\rmdefault}{\mddefault}{\updefault}{\color[rgb]{0,0,0}$u_{e_{3}}$}%
}}}}
\put(2786,781){\makebox(0,0)[lb]{\smash{{\SetFigFont{12}{14.4}{\rmdefault}{\mddefault}{\updefault}{\color[rgb]{0,0,0}$u^{1}_{3}$}%
}}}}
\put(3671, 31){\makebox(0,0)[lb]{\smash{{\SetFigFont{12}{14.4}{\rmdefault}{\mddefault}{\updefault}{\color[rgb]{0,0,0}$u^{2}_{3}$}%
}}}}
\put(2455,658){\makebox(0,0)[lb]{\smash{{\SetFigFont{12}{14.4}{\rmdefault}{\mddefault}{\updefault}{\color[rgb]{0,0,0}$u^{1}_{2}$}%
}}}}
\end{picture}%

%% file: anlexmpl.pdf_t
\begin{picture}(0,0)%
\includegraphics{anlexmpl.pdf}%
\end{picture}%
\setlength{\unitlength}{3947sp}%
\begingroup\makeatletter\ifx\SetFigFont\undefined%
\gdef\SetFigFont#1#2#3#4#5{%
  \reset@font\fontsize{#1}{#2pt}%
  \fontfamily{#3}\fontseries{#4}\fontshape{#5}%
  \selectfont}%
\fi\endgroup%
\begin{picture}(5416,1612)(218,-744)
\put(4564,685){\makebox(0,0)[lb]{\smash{{\SetFigFont{12}{14.4}{\rmdefault}{\mddefault}{\updefault}{\color[rgb]{0,0,0}$x^{j-1}_{P_{i}}$}%
}}}}
\put(503,-301){\makebox(0,0)[lb]{\smash{{\SetFigFont{12}{14.4}{\rmdefault}{\mddefault}{\updefault}{\color[rgb]{0,0,0}$y^{j}_{P_{i}}$}%
}}}}
\put(3354,-301){\makebox(0,0)[lb]{\smash{{\SetFigFont{12}{14.4}{\rmdefault}{\mddefault}{\updefault}{\color[rgb]{0,0,0}$y^{j}_{P_{i}}$}%
}}}}
\put(1710,685){\makebox(0,0)[lb]{\smash{{\SetFigFont{12}{14.4}{\rmdefault}{\mddefault}{\updefault}{\color[rgb]{0,0,0}$x^{j-1}_{P_{i}}$}%
}}}}
\end{picture}%

%% file: imppths.pdf_t
\begin{picture}(0,0)%
\includegraphics{imppths.pdf}%
\end{picture}%
\setlength{\unitlength}{3947sp}%
\begingroup\makeatletter\ifx\SetFigFont\undefined%
\gdef\SetFigFont#1#2#3#4#5{%
  \reset@font\fontsize{#1}{#2pt}%
  \fontfamily{#3}\fontseries{#4}\fontshape{#5}%
  \selectfont}%
\fi\endgroup%
\begin{picture}(2004,2047)(-449,-914)
\put(508,632){\makebox(0,0)[lb]{\smash{{\SetFigFont{12}{14.4}{\rmdefault}{\mddefault}{\updefault}{\color[rgb]{0,0,0}$T_{k}$}%
}}}}
\put(1118,411){\makebox(0,0)[lb]{\smash{{\SetFigFont{12}{14.4}{\rmdefault}{\mddefault}{\updefault}{\color[rgb]{0,0,0}$T_{2}$}%
}}}}
\put(1540,-247){\makebox(0,0)[lb]{\smash{{\SetFigFont{12}{14.4}{\rmdefault}{\mddefault}{\updefault}{\color[rgb]{0,0,0}$x_{i_{0}}^{j_{0}}$}%
}}}}
\put(1515,399){\makebox(0,0)[lb]{\smash{{\SetFigFont{12}{14.4}{\rmdefault}{\mddefault}{\updefault}{\color[rgb]{0,0,0}$x_{3}^{j_{0}}$}%
}}}}
\put(332,919){\makebox(0,0)[lb]{\smash{{\SetFigFont{12}{14.4}{\rmdefault}{\mddefault}{\updefault}{\color[rgb]{0,0,0}$x_{k}^{j_{0}}$}%
}}}}
\put(763,950){\makebox(0,0)[lb]{\smash{{\SetFigFont{12}{14.4}{\rmdefault}{\mddefault}{\updefault}{\color[rgb]{0,0,0}$x_{1}^{j_{0}}$}%
}}}}
\put(-253,509){\makebox(0,0)[lb]{\smash{{\SetFigFont{12}{14.4}{\rmdefault}{\mddefault}{\updefault}{\color[rgb]{0,0,0}$L_{j_{0}}$}%
}}}}
\put(841,615){\makebox(0,0)[lb]{\smash{{\SetFigFont{12}{14.4}{\rmdefault}{\mddefault}{\updefault}{\color[rgb]{0,0,0}$T_{1}$}%
}}}}
\put(1291,715){\makebox(0,0)[lb]{\smash{{\SetFigFont{12}{14.4}{\rmdefault}{\mddefault}{\updefault}{\color[rgb]{0,0,0}$x_{2}^{j_{0}}$}%
}}}}
\put(189,-836){\makebox(0,0)[lb]{\smash{{\SetFigFont{12}{14.4}{\rmdefault}{\mddefault}{\updefault}{\color[rgb]{0,0,0}$x_{i_{0}+1}^{j_{0}}$}%
}}}}
\put(-434,-157){\makebox(0,0)[lb]{\smash{{\SetFigFont{12}{14.4}{\rmdefault}{\mddefault}{\updefault}{\color[rgb]{0,0,0}$x_{i_{0}+2}^{j_{0}}$}%
}}}}
\put(750,-470){\makebox(0,0)[lb]{\smash{{\SetFigFont{12}{14.4}{\rmdefault}{\mddefault}{\updefault}{\color[rgb]{0,0,0}$T_{i_{0}}$}%
}}}}
\put(165,-340){\makebox(0,0)[lb]{\smash{{\SetFigFont{12}{14.4}{\rmdefault}{\mddefault}{\updefault}{\color[rgb]{0,0,0}$T_{i_{0}+1}$}%
}}}}
\end{picture}%

%% file: pthsexmpl.pdf_t
\begin{picture}(0,0)%
\includegraphics{pthsexmpl.pdf}%
\end{picture}%
\setlength{\unitlength}{3947sp}%
\begingroup\makeatletter\ifx\SetFigFont\undefined%
\gdef\SetFigFont#1#2#3#4#5{%
  \reset@font\fontsize{#1}{#2pt}%
  \fontfamily{#3}\fontseries{#4}\fontshape{#5}%
  \selectfont}%
\fi\endgroup%
\begin{picture}(12772,6889)(-2635,-4497)
\put(-1874,-2311){\makebox(0,0)[lb]{\smash{{\SetFigFont{12}{14.4}{\rmdefault}{\mddefault}{\updefault}{\color[rgb]{0,0,0}$L_{2}'=L_{j_{0}-2}$}%
}}}}
\put(976, 89){\makebox(0,0)[lb]{\smash{{\SetFigFont{12}{14.4}{\rmdefault}{\mddefault}{\updefault}{\color[rgb]{0,0,0}$u_{1}$}%
}}}}
\put(2176, 89){\makebox(0,0)[lb]{\smash{{\SetFigFont{12}{14.4}{\rmdefault}{\mddefault}{\updefault}{\color[rgb]{0,0,0}$u_{2}$}%
}}}}
\put(3376, 89){\makebox(0,0)[lb]{\smash{{\SetFigFont{12}{14.4}{\rmdefault}{\mddefault}{\updefault}{\color[rgb]{0,0,0}$u_{3}$}%
}}}}
\put(4576, 89){\makebox(0,0)[lb]{\smash{{\SetFigFont{12}{14.4}{\rmdefault}{\mddefault}{\updefault}{\color[rgb]{0,0,0}$u_{4}$}%
}}}}
\put(5776, 89){\makebox(0,0)[lb]{\smash{{\SetFigFont{12}{14.4}{\rmdefault}{\mddefault}{\updefault}{\color[rgb]{0,0,0}$u_{5}$}%
}}}}
\put(751,-1111){\makebox(0,0)[lb]{\smash{{\SetFigFont{12}{14.4}{\rmdefault}{\mddefault}{\updefault}{\color[rgb]{0,0,0}$L_{j_{0}}$}%
}}}}
\put(676,1739){\makebox(0,0)[lb]{\smash{{\SetFigFont{12}{14.4}{\rmdefault}{\mddefault}{\updefault}{\color[rgb]{0,0,0}$R_{1}$}%
}}}}
\put(1876,1739){\makebox(0,0)[lb]{\smash{{\SetFigFont{12}{14.4}{\rmdefault}{\mddefault}{\updefault}{\color[rgb]{0,0,0}$R_{2}$}%
}}}}
\put(3076,1739){\makebox(0,0)[lb]{\smash{{\SetFigFont{12}{14.4}{\rmdefault}{\mddefault}{\updefault}{\color[rgb]{0,0,0}$R_{3}$}%
}}}}
\put(4276,1739){\makebox(0,0)[lb]{\smash{{\SetFigFont{12}{14.4}{\rmdefault}{\mddefault}{\updefault}{\color[rgb]{0,0,0}$R_{4}$}%
}}}}
\put(5476,1739){\makebox(0,0)[lb]{\smash{{\SetFigFont{12}{14.4}{\rmdefault}{\mddefault}{\updefault}{\color[rgb]{0,0,0}$R_{5}$}%
}}}}
\put(-674,-2311){\makebox(0,0)[lb]{\smash{{\SetFigFont{12}{14.4}{\rmdefault}{\mddefault}{\updefault}{\color[rgb]{0,0,0}$L_{1}'=L_{j_{0}-1}$}%
}}}}
\end{picture}%

%% file: pthsexmpl2.pdf_t
\begin{picture}(0,0)%
\includegraphics{pthsexmpl2.pdf}%
\end{picture}%
\setlength{\unitlength}{3947sp}%
\begingroup\makeatletter\ifx\SetFigFont\undefined%
\gdef\SetFigFont#1#2#3#4#5{%
  \reset@font\fontsize{#1}{#2pt}%
  \fontfamily{#3}\fontseries{#4}\fontshape{#5}%
  \selectfont}%
\fi\endgroup%
\begin{picture}(12772,6889)(-2635,-4497)
\put(6676,-3511){\makebox(0,0)[lb]{\smash{{\SetFigFont{12}{14.4}{\rmdefault}{\mddefault}{\updefault}{\color[rgb]{0,0,0}$x_{4}^{j_{0}-1}$}%
}}}}
\put(976, 89){\makebox(0,0)[lb]{\smash{{\SetFigFont{12}{14.4}{\rmdefault}{\mddefault}{\updefault}{\color[rgb]{0,0,0}$u_{1}$}%
}}}}
\put(2176, 89){\makebox(0,0)[lb]{\smash{{\SetFigFont{12}{14.4}{\rmdefault}{\mddefault}{\updefault}{\color[rgb]{0,0,0}$u_{2}$}%
}}}}
\put(3376, 89){\makebox(0,0)[lb]{\smash{{\SetFigFont{12}{14.4}{\rmdefault}{\mddefault}{\updefault}{\color[rgb]{0,0,0}$u_{3}$}%
}}}}
\put(4576, 89){\makebox(0,0)[lb]{\smash{{\SetFigFont{12}{14.4}{\rmdefault}{\mddefault}{\updefault}{\color[rgb]{0,0,0}$u_{4}$}%
}}}}
\put(5776, 89){\makebox(0,0)[lb]{\smash{{\SetFigFont{12}{14.4}{\rmdefault}{\mddefault}{\updefault}{\color[rgb]{0,0,0}$u_{5}$}%
}}}}
\put(751,-1111){\makebox(0,0)[lb]{\smash{{\SetFigFont{12}{14.4}{\rmdefault}{\mddefault}{\updefault}{\color[rgb]{0,0,0}$L_{j_{0}}$}%
}}}}
\put(676,1739){\makebox(0,0)[lb]{\smash{{\SetFigFont{12}{14.4}{\rmdefault}{\mddefault}{\updefault}{\color[rgb]{0,0,0}$R_{1}$}%
}}}}
\put(1876,1739){\makebox(0,0)[lb]{\smash{{\SetFigFont{12}{14.4}{\rmdefault}{\mddefault}{\updefault}{\color[rgb]{0,0,0}$R_{2}$}%
}}}}
\put(3076,1739){\makebox(0,0)[lb]{\smash{{\SetFigFont{12}{14.4}{\rmdefault}{\mddefault}{\updefault}{\color[rgb]{0,0,0}$R_{3}$}%
}}}}
\put(4276,1739){\makebox(0,0)[lb]{\smash{{\SetFigFont{12}{14.4}{\rmdefault}{\mddefault}{\updefault}{\color[rgb]{0,0,0}$R_{4}$}%
}}}}
\put(5476,1739){\makebox(0,0)[lb]{\smash{{\SetFigFont{12}{14.4}{\rmdefault}{\mddefault}{\updefault}{\color[rgb]{0,0,0}$R_{5}$}%
}}}}
\put(-674,-2311){\makebox(0,0)[lb]{\smash{{\SetFigFont{12}{14.4}{\rmdefault}{\mddefault}{\updefault}{\color[rgb]{0,0,0}$L_{1}'=L_{j_{0}-1}$}%
}}}}
\put(-1874,-2311){\makebox(0,0)[lb]{\smash{{\SetFigFont{12}{14.4}{\rmdefault}{\mddefault}{\updefault}{\color[rgb]{0,0,0}$L_{2}'=L_{j_{0}-2}$}%
}}}}
\put(3076,-2461){\makebox(0,0)[lb]{\smash{{\SetFigFont{12}{14.4}{\rmdefault}{\mddefault}{\updefault}{\color[rgb]{0,0,0}$x_{1}^{j_{0}}$}%
}}}}
\put(1576,-2461){\makebox(0,0)[lb]{\smash{{\SetFigFont{12}{14.4}{\rmdefault}{\mddefault}{\updefault}{\color[rgb]{0,0,0}$x_{2}^{j_{0}}$}%
}}}}
\put(976,-2461){\makebox(0,0)[lb]{\smash{{\SetFigFont{12}{14.4}{\rmdefault}{\mddefault}{\updefault}{\color[rgb]{0,0,0}$x_{3}^{j_{0}}$}%
}}}}
\put(5776,-2461){\makebox(0,0)[lb]{\smash{{\SetFigFont{12}{14.4}{\rmdefault}{\mddefault}{\updefault}{\color[rgb]{0,0,0}$x_{4}^{j_{0}}$}%
}}}}
\put(4576,-2461){\makebox(0,0)[lb]{\smash{{\SetFigFont{12}{14.4}{\rmdefault}{\mddefault}{\updefault}{\color[rgb]{0,0,0}$x_{5}^{j_{0}}$}%
}}}}
\put(5776,-3661){\makebox(0,0)[lb]{\smash{{\SetFigFont{12}{14.4}{\rmdefault}{\mddefault}{\updefault}{\color[rgb]{0,0,0}$x_{5}^{j_{0}-2}$}%
}}}}
\put(3001,-4111){\makebox(0,0)[lb]{\smash{{\SetFigFont{12}{14.4}{\rmdefault}{\mddefault}{\updefault}{\color[rgb]{0,0,0}$x_{1}^{j_{0}-2}$}%
}}}}
\put(976,-3511){\makebox(0,0)[lb]{\smash{{\SetFigFont{12}{14.4}{\rmdefault}{\mddefault}{\updefault}{\color[rgb]{0,0,0}$x_{2}^{j_{0}-1}$}%
}}}}
\end{picture}%

%% file: surfim-arxiv.bbl
\begin{thebibliography}{10}

\bibitem{Abu-KhzamL03}
F.~N. Abu-Khzam and M.~A. Langston.
\newblock Graph coloring and the immersion order.
\newblock In T.~Warnow and B.~Zhu, editors, {\em COCOON}, volume 2697 of {\em
  Lecture Notes in Computer Science}, pages 394--403. Springer, 2003.

\bibitem{BelmonteHKPT12}
R.~Belmonte, P.~van~'t Hof, M.~Kaminski, D.~Paulusma, and D.~M. Thilikos.
\newblock Characterizing graphs of small carving-width.
\newblock In G.~Lin, editor, {\em COCOA}, volume 7402 of {\em Lecture Notes in
  Computer Science}, pages 360--370. Springer, 2012.

\bibitem{BiedlK97}
T.~C. Biedl and M.~Kaufmann.
\newblock Area-efficient static and incremental graph drawings.
\newblock In R.~E. Burkard and G.~J. Woeginger, editors, {\em ESA}, volume 1284
  of {\em Lecture Notes in Computer Science}, pages 37--52. Springer, 1997.

\bibitem{DemaineFHT05sube}
E.~D. Demaine, F.~V. Fomin, M.~Hajiaghayi, and D.~M. Thilikos.
\newblock Subexponential parameterized algorithms on bounded-genus graphs and
  {$H$}-minor-free graphs.
\newblock {\em J.~ACM}, 52(6):866--893, 2005.

\bibitem{2011arXiv1101.2630D}
M.~{DeVos}, Z.~{Dvo{\v r}{\'a}k}, J.~{Fox}, J.~{McDonald}, B.~{Mohar}, and
  D.~{Scheide}.
\newblock {Minimum degree condition forcing complete graph immersion}.
\newblock {\em ArXiv e-prints}, Jan. 2011.

\bibitem{1213.05137}
M.~DeVos, K.-I. Kawarabayashi, B.~Mohar, and H.~Okamura.
\newblock Immersing small complete graphs.
\newblock {\em Ars Math. Contemp.}, 3(2):139--146, 2010.

\bibitem{FerraraGTW08}
M.~Ferrara, R.~J. Gould, G.~Tansey, and T.~Whalen.
\newblock On {\it h}-immersions.
\newblock {\em Journal of Graph Theory}, 57(3):245--254, 2008.

\bibitem{FominGT11}
F.~V. Fomin, P.~A. Golovach, and D.~M. Thilikos.
\newblock Contraction obstructions for treewidth.
\newblock {\em J. Comb. Theory, Ser. B}, 101(5):302--314, 2011.

\bibitem{kurim}
A.~C. Giannopoulou, M.~Kaminski, and D.~M. Thilikos.
\newblock Forbidding kuratowski graphs as immersions.
\newblock {\em CoRR}, abs/1207.5329, 2012.

\bibitem{GiannopoulouSZ12}
A.~C. Giannopoulou, I.~Salem, and D.~Zoros.
\newblock Effective computation of immersion obstructions for unions of graph
  classes.
\newblock In F.~V. Fomin and P.~Kaski, editors, {\em SWAT}, volume 7357 of {\em
  Lecture Notes in Computer Science}, pages 165--176. Springer, 2012.

\bibitem{GiannopoulouT11}
A.~C. Giannopoulou and D.~M. Thilikos.
\newblock Optimizing the graph minors weak structure theorem.
\newblock {\em CoRR}, abs/1102.5762, 2011.

\bibitem{GroheKMW11}
M.~Grohe, K.~ichi Kawarabayashi, D.~Marx, and P.~Wollan.
\newblock {F}inding topological subgraphs is fixed-parameter tractable.
\newblock In {\em STOC}, pages 479--488, 2011.

\bibitem{KawarabayashiK12}
K.~ichi Kawarabayashi and Y.~Kobayashi.
\newblock List-coloring graphs without subdivisions and without immersions.
\newblock In Y.~Rabani, editor, {\em SODA}, pages 1425--1435. SIAM, 2012.

\bibitem{Kostochka84}
A.~V. Kostochka.
\newblock Lower bound of the hadwiger number of graphs by their average degree.
\newblock {\em Combinatorica}, 4(4):307--316, 1984.

\bibitem{Lescure1988325}
F.~Lescure and H.~Meyniel.
\newblock On a problem upon configurations contained in graphs with given
  chromatic number.
\newblock In L.~D. Andersen, I.~T. Jakobsen, C.~Thomassen, B.~Toft, and P.~D.
  Vestergaard, editors, {\em Graph Theory in Memory of G.A. Dirac}, volume~41
  of {\em Annals of Discrete Mathematics}, pages 325 -- 331. Elsevier, 1988.

\bibitem{MoharT01}
B.~Mohar and C.~Thomassen.
\newblock {\em Graphs on surfaces}.
\newblock Johns Hopkins University Press, Baltimore, MD, 2001.

\bibitem{RobertsonST94}
N.~Robertson, P.~Seymour, and R.~Thomas.
\newblock Quickly excluding a planar graph.
\newblock {\em J. Combin. Theory Ser. B}, 62(2):323--348, 1994.

\bibitem{RobertsonS-GMXIII}
N.~Robertson and P.~D. Seymour.
\newblock Graph minors. {XIII}. {T}he disjoint paths problem.
\newblock {\em J. Combin. Theory Ser. B}, 63(1):65--110, 1995.

\bibitem{RobertsonS-XVI}
N.~Robertson and P.~D. Seymour.
\newblock Graph minors. {X}{V}{I}. {E}xcluding a non-planar graph.
\newblock {\em J. Comb. Theory, Ser. B}, 89(1):43--76, 2003.

\bibitem{RobertsonS10}
N.~Robertson and P.~D. Seymour.
\newblock {G}raph {m}inors {X}{X}{I}{I}{I}. {N}ash-{W}illiams' immersion
  conjecture.
\newblock {\em J. Comb. Theory, Ser. B}, 100(2):181--205, 2010.

\bibitem{Thomason01}
A.~Thomason.
\newblock The extremal function for complete minors.
\newblock {\em J. Combin. Theory Ser. B}, 81(2):318--338, 2001.

\bibitem{abs-1302-3867}
P.~Wollan.
\newblock The structure of graphs not admitting a fixed immersion.
\newblock {\em CoRR}, abs/1302.3867, 2013.

\end{thebibliography}
